\documentclass[12pt]{article}
\usepackage{amssymb}
\usepackage{amsmath}
\usepackage{array}
\usepackage[russian]{babel}
\usepackage[pdftex]{graphics}

\usepackage[X2,T2A]{fontenc}
\usepackage[cp1251]{inputenc}

\DeclareMathOperator{\vep}{\varepsilon}

\newcommand{\Z}{\scriptstyle}

\newcommand{\prsum}{\mathop{{\sum}'}}

\multlinegap=0em \DeclareFontFamily{T1}{msb}{}
\DeclareFontShape{T1}{msb}{m}{ol}{<5> <6> <7> <8> <9> gen * msbm
<10> <10.95> <12> <14.4> <17.28> <20.74> <24.88> msbm10}{}
\DeclareSymbolFont{AMSb}{T1}{msb}{m}{ol} \multlinegap=0em

\renewcommand{\S}{\mathhexbox278}
\renewcommand{\le}{\operatorname{\leqslant}}
\renewcommand{\ge}{\operatorname{\geqslant}}

\begin{document}

\begin{flushleft}
MSC 11L05, 11N37
\end{flushleft}

\begin{center}
{\rmfamily\bfseries\normalsize Short Kloosterman sums to powerful modulus}
\end{center}

\begin{center}
{\rmfamily\bfseries\normalsize M.A.~Korolev\footnote{The work is supported by the Russian Science Foundation under grant 14-11-00433 and performed in Steklov Mathematical Institure of Russian Academy of Sciences}}
\end{center}

\vspace{0.5cm}

\fontsize{11}{12pt}\selectfont

\textbf{Abstract.} We obtain the estimate of incomplete Kloosterman sum to powerful modulus $q$. The length $N$ of the sum
lies in the interval $e^{c(\log{q})^{2/3}}\le N\le \sqrt{q}$.

\vspace{0.2cm}


\fontsize{12}{15pt}\selectfont

\vspace{0.5cm}

The aim of this paper is to estimate a short Kloosterman sum
\begin{equation}\label{lab_01}
S\,=\,S_{q}(N;a,b,c)\,=\,\prsum\limits_{c<n\le c+N}e_{q}(an^{*}+bn),
\end{equation}
to powerful modulus $q$. Here $N,a,b,c$ are integers, $1<N\le\sqrt{q}$, $(a,q)=1$, the prime sign in the sum means the summation over $(n,q)=1$, $nn^{*}\equiv 1 \pmod{q}$ and $e_{q}(v) = e^{2\pi iv/q}$. The number $q$ is called powerful if its kernel $d = \prod_{p|q}p$ is small relative to $q$ in the logarithmic scale. The simplest case of such numbers are $q = p^{n}$ where $p$ is fixed prime and $n\to +\infty$.

For many times, A.A.~Karatsuba \cite{Karatsuba_1994}-\cite{Karatsuba_1997} pointed out to the possibility of estimating of such sums with $q = p^{n}$ by method of A.G.~Postnikov \cite{Postnikov_1955}, \cite{Postnikov_1956}. In this paper, we prove the following statement.

\vspace{0.3cm}

\textsc{Theorem 1.} \emph{Suppose that $q\ge q_{0}$ is sufficiently large, $d = \prod_{p|q}p$ is the kernel of $q$, $\gamma_{1} = 900$, $\gamma = 160^{-4}$ and let $\max{\bigl(d^{15},e^{\gamma_{1}(\ln{q})^{2/3}}\bigr)}\le N\le \sqrt{q}$. Then, for any $a,b,c$ such that $(a,q)=1$, the following estimate holds:}
\begin{equation}\label{lab_02}
|S_{q}(N;a,b,c)|\,<\,N\exp{\biggl(-\,\gamma\,\frac{(\ln{N})^{3}}{(\ln{q})^{2}}\biggr)}.
\end{equation}

\vspace{0.3cm}

The proof does not use any new ideas and is based only on the technic of papers \cite{Karatsuba_1964}, \cite{Chubarikov_1973} and \cite{Iwaniec_1974}.

We need the following auxilliary assertions.

\vspace{0.3cm}

\textsc{Lemma 1.} \emph{Suppose that} $0<\vep<1$, $m = \bigl[2\vep^{-1}\bigr]$, $q=p_{1}^{\alpha_{1}}\ldots p_{s}^{\alpha_{s}}$,
$d = \prod_{p|q}p$, \emph{and let $q_{\vep} = dp_{1}^{\beta_{1}}\ldots p_{s}^{\beta_{s}}$, where $\beta_{r}=\bigl[\vep\alpha_{r}\bigr]$. Then for any $z$ we have}
\begin{equation}\label{lab_03}
(1+zq_{\vep})^{*}\,\equiv\,1-zq_{\vep}+(zq_{\vep})^{2}-\ldots + (-1)^{m}(zq_{\vep})^{m}\pmod{q}.
\end{equation}

\textsc{Remark.} This assertion is an analogue of Postnikov's formula for $\text{ind}\,(1+pz)$ modulo $q = p^{n}$, $p\ge 3$ (see \cite{Postnikov_1955}, \cite{Postnikov_1956}). The idea of introducing the factor $q_{\vep}$ belongs to H.~Iwaniec \cite{Iwaniec_1974}.

\vspace{0.3cm}

\textsc{Proof.} It is sufficient to prove that $q_{\vep}^{m+1}$ is divisible by $q$. This fact follows from easy-to-check inequalities $(m+1)(\beta_{r}+1)\ge\alpha_{r}$, $r = 1,2,\ldots,s$.

\vspace{0.3cm}

\textsc{Lemma 2.} \emph{Suppose that $P\ge 1$ and $\alpha$ are any real numbers. Then}
\[
\biggl|\sum\limits_{1\le n\le P}e(\alpha n)\biggr|\,\le\,\min{\bigl(P,\|\alpha\|^{-1}\bigr)},
\]
\emph{where} $e(z) = e^{2\pi iz}$ \emph{and} $\|\alpha\| = \min{\bigl(\{\alpha\},1-\{\alpha\}\bigr)}$.

\vspace{0.3cm}

\textsc{Lemma 3.} \emph{Suppose that}
\[
\alpha\,=\,\frac{A}{Q}\,+\,\frac{\theta}{Q^{2\mathstrut}},\quad (A,Q)=1,\quad Q\ge 1,\quad |\theta|\le 1.
\]
\emph{and let $\beta$, $U>0$, $P\ge 1$ be any real numbers. Then the following inequality holds:}
\[
\sum\limits_{1\le n\le P}\min{\bigl(U,\|\alpha n+\beta\|^{-1}\bigr)}\,\le\,6\biggl(\frac{P}{Q}+1\biggr)(U+Q\log{Q}).
\]

For the proofs of lemmas 2,3, see, for example, \cite[Ch. VI, \S 2]{Karatsuba_1993}.

Denote by $J_{k,m}(P;\lambda_{1},\ldots,\lambda_{m})$ the number of solutions of the following system:
\begin{equation*}
\begin{cases}
& x_{1}+\ldots+x_{k}\,=\,x_{k+1} + \ldots + x_{2k} + \lambda_{1}\\
& \cdots\\
& x_{1}^{m}+\ldots+x_{k}^{m}\,=\,x_{k+1}^{m} + \ldots + x_{2k}^{m} + \lambda_{m}
\end{cases}
\end{equation*}
with integers variables $1\le x_{1},\ldots, x_{2k}\le P$ and let $J_{k,m}(P)=J_{k,m}(P;0,\ldots,0)$. Obviously,
$J_{k,m}(P;\lambda_{1},\ldots,\lambda_{m})\le J_{k,m}(P)$ for any $\lambda_{1},\ldots,\lambda_{m}$.

\vspace{0.3cm}

\textsc{Lemma 4 (Vinogradov's mean value theorem).} \emph{Suppose that $\tau\ge 1$, $k\ge m\tau$, $P\ge 1$ are integers. Then}
$J_{k,m}(P)\,\le\,D(m,\tau)P^{2k-\Delta(m,\tau)}$, \emph{where}
\[
D(m,\tau)\,=\,(m\tau)^{6m\tau}(2m)^{4m(m+1)\tau},\quad \Delta(m,\tau)\,=\,\frac{1}{2}\,m(m+1)\biggl(1\,-\,\biggl(1\,-\,\frac{1}{m}\biggr)^{\!\tau}\,\biggr).
\]

For the proof of this version of Vinogradov's mean value theorem, see \cite[Ch. VI, \S 1]{Karatsuba_1993}.

\vspace{0.3cm}

\textsc{Proof of Theorem 1.} Shifting the interval of summation to at most $d\le N^{1/15}$ terms, we obtain the inequality
$|S|\le |S_{1}|+d$, where the sum $S_{1}$ has the same type as the initial sum $S$ and satisfies the additional condition $c\equiv 0\pmod{d}$.

Further, suppose that $h = [N^{1/4}]+1$, $1\le x,y\le h$ and define $q_{\vep}$ as in lemma 1 (the precise value of $\vep$ will be chosen later). Then
\begin{multline*}
S_{1}\,=\,\prsum\limits_{c<n+q_{\vep}xy\le c+N}e_{q}\bigl(a(n+q_{\vep}xy)^{*}+b(n+q_{\vep}xy)\bigr)\,=\\
=\,\prsum\limits_{1\le n\le N}e_{q}\bigl(a(n+c+q_{\vep}xy)^{*}+b(n+c+q_{\vep}xy)\bigr)\,+\,2\theta q_{\vep}xy,\quad |\theta|\le 1.
\end{multline*}
Summing over $1\le x,y\le h$ we get
\[
|S_{1}|\,\le\,h^{-2}\prsum\limits_{1\le n\le N}|W|\,+\,h^{2}q_{\vep},\
\]
where
\[
W\,=\,W(n)\,=\,\sum\limits_{x,y=1}^{h}e_{q}\bigl(a(n+c+q_{\vep}xy)^{*}+bq_{\vep}xy\bigr).
\]
Setting $v \equiv (n+c)^{*}\pmod{q}$ for brevity and using lemma 1, we obtain
\begin{multline*}
|W|\,=\,\biggl|\sum\limits_{x,y=1}^{h}e_{q}\bigl(a_{1}xy+a_{2}(xy)^{2}+\ldots+a_{m}(xy)^{m}\bigr)\biggr|\,=\\
=\,\biggl|\sum\limits_{x,y=1}^{h}e\bigl(\alpha_{1}xy+\alpha_{2}(xy)^{2}+\ldots+\alpha_{m}(xy)^{m}\bigr)\biggr|,
\end{multline*}
where $e(z) = e^{2\pi iz}$, $m\,=\,\bigl[2\vep^{-1}\bigr]$ and
\[
a_{1}\equiv q_{\vep}(b-av^{2})\pmod{q},\quad a_{r}\equiv (-1)^{r}av^{r+1}q_{\vep}^{r}\pmod{q},\quad r = 2,3,\ldots,m,\quad
\alpha_{r} = \frac{a_{r}}{q}.
\]
Taking $k = m\tau$ (where the value of $\tau$ will be chosen later) and applying H\"{o}lder's inequality to $W$ as in \cite[Ch. VI, \S 1]{Karatsuba_1993},
we find that
\begin{multline*}
|W|^{2k}\,\le\,h^{2k-1}\sum\limits_{x = 1}^{h}\biggl|\sum\limits_{y=1}^{h}e(\alpha_{1}xy+\ldots + \alpha_{m}x^{m}y^{m})\biggr|^{2k}\,=\\
=\,h^{2k-1}\sum\limits_{x = 1}^{h}\sum\limits_{\lambda_{1},\ldots,\lambda_{m}}J_{k,m}(h;\lambda_{1},\ldots,\lambda_{m})e(\alpha_{1}\lambda_{1}x+\ldots + \alpha_{m}\lambda_{m}x^{m})\,\le\\
\le\,h^{2k-1}\sum\limits_{\lambda_{1},\ldots,\lambda_{m}}J_{k,m}(h;\lambda_{1},\ldots,\lambda_{m})\biggl|\sum\limits_{x=1}^{h}e(\alpha_{1}\lambda_{1}x+\ldots + \alpha_{m}\lambda_{m}x^{m})\biggr|,
\end{multline*}
where $\lambda_{r}$ runs through some set of values, $|\lambda_{r}|<\Lambda_{r}$, $\Lambda_{r} = kh^{r}$. Obviously,
\[
\sum\limits_{\lambda_{1},\ldots,\lambda_{m}}J_{k,m}(h;\lambda_{1},\ldots,\lambda_{m})\,\le\,h^{2k}.
\]
Using H\"{o}lder's inequality again, we obtain:
\begin{multline*}
|W|^{4k^{2}}\,\le\,h^{2k(2k-1)}\biggl(\sum\limits_{\lambda_{1},\ldots,\lambda_{m}}J_{k,m}(h;\lambda_{1},\ldots,\lambda_{m})\biggr)^{\!2k-1}\times\\
\times \sum\limits_{\lambda_{1},\ldots,\lambda_{m}}J_{k,m}(h;\lambda_{1},\ldots,\lambda_{m})\biggl|\sum\limits_{x=1}^{h}e(\alpha_{1}\lambda_{1}x+\ldots + \alpha_{m}\lambda_{m}x^{m})\biggr|^{2k}\,\le\\
\le\,h^{4k(2k-1)}J_{k,m}(h)\sum\limits_{\lambda_{1},\ldots,\lambda_{m}}\biggl|\sum\limits_{x=1}^{h}e(\alpha_{1}\lambda_{1}x+\ldots + \alpha_{m}\lambda_{m}x^{m})\biggr|^{2k},
\end{multline*}
where $\lambda_{r}$ in the last sum runs through the entire interval $|\lambda_{r}|<\Lambda_{r}$ ($r = 1,\ldots,m$). Further,
\begin{multline*}
|W|^{4k^{2}}\,\le\\
\le\,h^{4k(2k-1)}J_{k,m}(h)\sum\limits_{\lambda_{1},\ldots,\lambda_{m}}\sum\limits_{\mu_{1},\ldots,\,\mu_{m}}J_{k,m}(h;\mu_{1},\ldots,\mu_{m})
e(\alpha_{1}\lambda_{1}\mu_{1}+\ldots + \alpha_{m}\lambda_{m}\mu_{m})\,\le\\
\le\,h^{4k(2k-1)}J_{k,m}(h)\sum\limits_{\mu_{1},\ldots,\,\mu_{m}}J_{k,m}(h;\mu_{1},\ldots,\mu_{m})\biggl|\sum\limits_{\lambda_{1}}e(\alpha_{1}\mu_{1}\lambda_{1})\biggr|
\cdots \biggl|\sum\limits_{\lambda_{m}}e(\alpha_{m}\mu_{m}\lambda_{m})\biggr|.
\end{multline*}
By lemma 2,
\[
|W|^{4k^{2}}\,\le\,h^{4k(2k-1)}J_{k,m}^{\,2}(h)V_{1}\cdots V_{m},
\]
where
\[
V_{r}\,=\,\sum\limits_{|\mu_{r}|<\Lambda_{r}}\min{\bigl(2\Lambda_{r},\|\alpha_{r}\mu_{r}\|^{-1}\bigr)},\quad r = 1,2,\ldots, m.
\]
Obviously, the trivial bound for $V_{r}$ is $(2\Lambda_{r})^{2}$. At the same time, representing $\alpha_{r}$ as incontractible fraction of the form
\begin{equation}\label{lab_04}
\frac{A_{r}}{Q_{r}},\quad Q_{r}\,=\,p_{1}^{\gamma_{1}}\ldots p_{s}^{\gamma_{s}},\quad \gamma_{\ell}\,=\,\max{(0,\alpha_{\ell}-r\beta_{\ell})},
\end{equation}
we obtain:
\begin{multline*}
V_{r}\,=\,6\biggl(\frac{2\Lambda_{r}}{Q_{r}}\,+\,1\biggr)(2\Lambda_{r}\,+\,Q_{r}\ln{Q_{r}})\,<\,(2\Lambda_{r})^{2}\delta_{r},\\
\delta_{r}\,=\,6(\log{q})\biggl(\frac{1}{Q_{r}}\,+\,\frac{1}{2\Lambda_{r}}\biggr)\biggl(1\,+\,\frac{Q_{r}}{2\Lambda_{r}}\biggr)\,=\,6(\log{q})
\biggl(\frac{1}{\sqrt{Q_{r}\mathstrut}}\,+\,\frac{\sqrt{Q_{r}\mathstrut}}{2\Lambda_{r}}\biggr)^{\!2}.
\end{multline*}
Hence, setting $\Delta_{r} = \min{(1,\delta_{r})}$, we get
\[
|W|^{4k^{2}}\,\le\,h^{4k(2k-1)}J_{k,m}^{\,2}(h)\Delta\prod\limits_{r=1}^{m}(2\Lambda_{r})^{2}\,=\,(2k)^{2m}h^{4k(2k-1)+m(m+1)}J_{k,m}^{\,2}(h)\Delta,
\]
where $\Delta\,=\,\Delta_{1}\cdots \Delta_{m}$. Using lemma 4, we obtain
\begin{multline*}
J_{k,m}^{\,2}(h)\,\le\,(m\tau)^{12m\tau}(2m)^{4m(m+1)\tau}h^{4k-m(m+1)\bigl(1-\bigl(1-\tfrac{\Z 1}{\Z m}\bigr)^{\!\tau}\bigr)}\Delta\,=\\
=\,k^{12k}(2m)^{4k(m+1)}h^{4k-m(m+1)\bigl(1-\bigl(1-\tfrac{\Z 1}{\Z m}\bigr)^{\!\tau}\bigr)}\Delta,\\
|W|^{4k^{2}}\,\le\,C(k,m)h^{8k^{2}+m(m+1)\bigl(1-\tfrac{\Z 1}{\Z m}\bigr)^{\!\tau}}\Delta,\quad C(k,m)\,=\,k^{12k}(2m)^{4k(m+1)}(2k)^{2m}.
\end{multline*}
Now let us choose
\[
\vep\,=\,c\,\frac{\log{N}}{\log{q}},\quad c\,=\,\frac{1}{7},\quad r_{j}\,=\,\bigl[c_{j}\vep^{-1}\bigr],\quad j=1,2,\quad
c_{1}\,=\,\frac{1}{3},\;c_{2}\,=\,\frac{2}{3},
\]
and note that $N = q^{7\vep}$, $r_{1}\ge 1$, $r_{2}-r_{1}\ge 1$. Suppose now that $r_{1}<r\le r_{2}$.
Then, by (\ref{lab_04}), we have $qq_{\vep}^{-r}\le Q_{r}\le qq_{\vep}^{-1}<q$ and hence
\[
\frac{1}{Q_{r}}\,\le\,\frac{q_{\vep}^{r}}{q}\,\le\,\frac{(dq_{\vep})^{r}}{q}\,\le\,\frac{(dq_{\vep})^{r_{2}}}{q}.
\]
Obviously, $q^{\,\vep r_{2}}\le q^{2/3}$, $d^{\,r_{2}}\le N^{r_{2}/15} = q^{7\vep r_{2}/15}\,\le\,q^{14/15}$,
so we have $Q_{r}\le q^{-1/45}$. Similarly,
\[
\frac{Q_{r}}{4\Lambda_{r}^{2\mathstrut}}\,<\,\frac{q}{h^{2r}}\,\le\,\frac{q}{h^{2(r_{1}+1)}}\le \frac{q}{N^{(r_{1}+1)/2}}\,\le\,q^{-1/6}.
\]
Hence,
\begin{multline*}
\delta_{r}\,\le\,6(\log{q})\bigl(q^{-1/90}\,+\,q^{-1/12}\bigr)^{2}\,<\,q^{-1/50},\\
\Delta\,\le\,\prod\limits_{r = r_{1}+1}^{r_{2}}\delta_{r}\,<\,q^{-(r_{2}-r_{1})/50}\,\le\,q^{-\vep^{-1}/200}\,=\,\exp{\biggl(-\,\frac{7}{200}\,\frac{(\log{q})^{2}}{\log{N}}\biggr)},\\
|W|^{4k^{2}}\,<\,C(k,m)h^{8k^{2}+m(m+1)\bigl(1-\tfrac{\Z 1}{\Z m}\bigr)^{\!\tau}}\exp{\biggl(-\,\frac{7}{200}\,\frac{(\log{q})^{2}}{\log{N}}\biggr)}.
\end{multline*}
Now we set $\tau = \kappa m$, $\kappa = 10$. Then we have
\begin{multline*}
m(m+1)\biggl(1\,-\,\frac{1}{m}\biggr)^{\!\tau}\,\le\,e^{-10}m(m+1)\,<\,\frac{1}{91}\,\frac{(\log{q})^{2}}{\log{N}},\\
h^{m(m+1)\bigl(1-\tfrac{\Z 1}{\Z m}\bigr)^{\!\tau}}\,\le\,\exp{\biggl(\frac{1}{360}\,\frac{(\log{q})^{2}}{\log{N}}\biggr)},\\
|W|^{4k^{2}}\,<\,C(k,m)h^{8k^{2}}\exp{\biggl(-\,\frac{1}{31.5}\,\frac{(\log{q})^{2}}{\log{N}}\biggr)}.
\end{multline*}
Since $0<\vep<\tfrac{1}{14}$ then $m\ge 28$ and
\[
C(k,m)^{1/(4k^{2})}\,=\,(10m^{2})^{\frac{\Z 1}{\Z 30m^{2\mathstrut}}}(20m^{2})^{\frac{\Z 1}{\Z 200m^{3\mathstrut}}}(2m)^{\frac{\Z m+1}{\Z 10m^{2\mathstrut}}}\,<\,1.02.
\]
Thus we have
\[
|W|\,<\,1.02h^{2}\exp{\biggl(-\,\frac{1}{126k^{2\mathstrut}}\,\frac{(\log{N})^{3}}{(\log{q})^{2\mathstrut}}\biggr)}\,<
\,1.02h^{2}\exp{\biggl(-\,\frac{1}{156^{4\mathstrut}}\,\frac{(\log{N})^{3}}{(\log{q})^{2\mathstrut}}\biggr)},
\]
and, finally,
\[
|S_{1}|\,<\,1.02N\exp{\biggl(-\,\frac{1}{156^{4\mathstrut}}\,\frac{(\log{N})^{3}}{(\log{q})^{2\mathstrut}}\biggr)},\quad
|S|\,\le\,|S_{1}|+d\,<\,N\exp{\biggl(-\,\frac{1}{160^{4\mathstrut}}\,\frac{(\log{N})^{3}}{(\log{q})^{2\mathstrut}}\biggr)}.
\]
Theorem is proved. $\Box$

\vspace{0.3cm}

The restriction $d^{15}\le N$ can be weaken slightly, but the price is the narrowing of the interval for $N$ and the loss of precision of the estimate. Namely, the following assertion is true.

\vspace{0.3cm}

\textsc{Theorem 2.} \emph{Suppose that $0<\delta<0.1$ is any fixed number, $q\ge q_{0}(\delta)$ and let}
\[
\gamma_{1}\,=\,1200\delta^{-2}\bigl(\log{(1/\delta)}\bigr)^{\!2/3},\quad \max{\bigl(d^{\,2+\delta},e^{\gamma_{1}(\log{q})^{2/3}}\bigr)}\le N\le q^{\delta/20}.
\]
\emph{Then the inequality} (\ref{lab_02}) \emph{holds with} $\gamma\,=\,201^{-4}\delta^{\,6}\bigl(\log{(1/\delta)}\bigr)^{\!2}$.

\vspace{0.3cm}

\textsc{Proof.} Setting $m = \bigl[2\vep^{-1}\bigr]$, $r_{j} = \bigl[c_{j}\vep^{-1}\bigr]$, $j=1,2$,
\[
\vep\,=\,c\,\frac{\log{N}}{\log{q}}\,\quad c\,=\,\frac{\delta}{5}\biggl(1\,-\,\frac{\delta}{15}\biggr),\quad c_{1}\,=\,\frac{2\delta}{5}\biggl(1\,-\,\frac{\delta}{20}\biggr),
\quad c_{2}\,=\,\frac{2\delta}{5}
\]
in the above calculations and noting that $q^{\vep} = N^{c}$, $d\le N^{1/(2+\delta)}$, we find
\[
q_{\vep}h^{2}\,\le\, 2\sqrt{N}dq^{\vep}\,\le\, 2N^{\xi},
\]
where
\[
\xi\,=\,\frac{1}{2}+\frac{1}{2+\delta}+c\,<\,\frac{1}{2}+\frac{1}{2}\biggl(1-\frac{\delta}{2}+\frac{\delta^{2}}{4}\biggr)+\frac{\delta}{5}
\,=\,1-\frac{\delta}{20}\biggl(1-\frac{5\delta}{2}\biggr)\,<\,1-\frac{\delta}{40}.
\]
Further, if $r_{1}<r\le r_{2}$ then
\begin{multline*}
\frac{1}{Q_{r}}\,\le\,\frac{(dq^{\vep})^{r_{2}}}{q}\,=\,q^{\eta},\\
\eta\,=\,-1+\vep r_{2}\biggl(1\,+\,\frac{1}{c(2+\delta)}\biggr)\,\le\,-1+c_{2}\biggl(1\,+\,\frac{1}{c(2+\delta)}\biggr)\,=\\
=\,\frac{2\delta}{5}\,+\,\frac{1-(1-\tfrac{\Z \delta}{\Z 15\mathstrut})(1+\tfrac{\Z \delta}{\Z 2\mathstrut})}{(1-\tfrac{\Z \delta}{\Z 15\mathstrut})(1+\tfrac{\Z \delta}{\Z 2\mathstrut})}\,=\,-\,\frac{\delta}{30}\,\frac{1-\tfrac{\Z 31\delta}{\Z 5\mathstrut}+\frac{\Z 2\delta^{2}}{\Z 5\mathstrut}}{(1-\tfrac{\Z \delta}{\Z 15\mathstrut})(1+\tfrac{\Z \delta}{\Z 2\mathstrut})}\,<\,-\,\frac{\delta}{150},
\end{multline*}
\begin{multline*}
\frac{Q_{r}}{4\Lambda_{r}^{2\mathstrut}}\,<\,\frac{q}{h^{2r}}\,\le\,\frac{q}{N^{(r_{1}+1)/2}}\,=\,q^{\vartheta},\\
\vartheta\,=\,1\,-\,\frac{\vep}{2c}\,(r_{1}+1)\,\le\,1\,-\,\frac{c_{1}}{2c}\,=\,-\,\frac{\delta}{60}\,\frac{1}{1-\tfrac{\Z \delta}{\Z 15\mathstrut}}\,<\,-\,\frac{\delta}{60},
\end{multline*}
so we have
\[
\delta_{r}\,<\,24(\log{q})q^{-\,\delta/150}\,<\,q^{-\,\delta/160},\quad \Delta\,<\,q^{-\delta (r_{2}-r_{1})/160}\,<\,
\exp{\biggl(-\,\frac{\delta^{2}}{3200}\,\frac{(\log{q})^{2}}{\log{N}}\biggr)}.
\]
Taking $\tau = \kappa m$, $\kappa = \bigl[4\log{(1/\delta)}\bigr]+14$ we find that
\[
h^{m(m+1)\bigl(1-\tfrac{\Z 1}{\Z m}\bigr)^{\!\tau}}\,\le\,h^{e^{-\kappa}m(m+1)}\,<\,\exp{\biggl(\frac{5e^{-\kappa}}{4c^{2\mathstrut}}\,\frac{(\log{q})^{2}}{\log{N}}\biggr)}
\,<\,\exp{\biggl(\frac{\delta^{2}}{6400}\,\frac{(\log{q})^{2}}{\log{N}}\biggr)}
\]
and hence
\[
|W|^{4k^{2}}\,<\,C(k,m)h^{8k^{2}}\exp{\biggl(-\,\frac{\delta^{2}}{6400}\,\frac{(\log{q})^{2}}{\log{N}}\biggr)},
\]
where $C(k,m)$ is defined above. Now it remains to check the inequalities
\[
\frac{\delta^{2}}{4k^{2\mathstrut}\cdot 6400}\,=\,\frac{\delta^{\,2}\kappa^{-2}c^{4}}{640^{2\mathstrut}}\,
\frac{(\log{q})^{4}}{(\log{N})^{4}}\,>\,\frac{\delta^{\,6}}{200^{4}}\bigl(\log{(1/\delta)}\bigr)^{-2}\,
\frac{(\log{q})^{4}}{(\log{N})^{4}}.
\]
Theorem 2 is proved. $\Box$

\renewcommand{\refname}{\normalsize{Bibliography}}

\textsc{Maxim A. Korolev}

Steklov Mathematical Institute of Russian Academy of Sciences

119991, Russia, Moscow, Gubkina str., 8

e\,-mail: \texttt{korolevma@mi.ras.ru}

\end{document}